\documentclass[12pt,reqno]{article}

\usepackage[usenames]{color}
\usepackage{amssymb}
\usepackage{amsmath}
\usepackage{amsthm}
\usepackage{amsfonts}
\usepackage{amscd}
\usepackage{graphicx}
\DeclareMathOperator{\lcm}{lcm}

\usepackage[colorlinks=true,
linkcolor=webgreen,
filecolor=webbrown,
citecolor=webgreen]{hyperref}

\definecolor{webgreen}{rgb}{0,.5,0}
\definecolor{webbrown}{rgb}{.6,0,0}

\usepackage{color}
\usepackage{fullpage}
\usepackage{float}

\usepackage{graphics}
\usepackage{latexsym}
\usepackage{epsf}
\usepackage{breakurl}

\def\lcm{\mathrm{lcm}}

\begin{document}

\begin{center}
\epsfxsize=4in
\end{center}

\theoremstyle{plain}
\newtheorem{theorem}{Theorem}
\newtheorem{corollary}[theorem]{Corollary}
\newtheorem{lemma}[theorem]{Lemma}
\newtheorem{proposition}[theorem]{Proposition}

\theoremstyle{definition}
\newtheorem{definition}[theorem]{Definition}
\newtheorem{example}[theorem]{Example}
\newtheorem{conjecture}[theorem]{Conjecture}

\theoremstyle{remark}
\newtheorem{remark}[theorem]{Remark}

\begin{center}
\vskip 1cm{\LARGE\bf  Extreme Covering Systems
}

\vskip 1cm
\large
Jack Dalton and Ognian Trifonov\\
University of South Carolina\\
Department of Mathematics \\
Columbia, SC 29208 \\
\href{mailto:jrdalton@math.sc.edu}{\tt jrdalton@math.sc.edu} \\
\href{mailto:trifonov@math.sc.edu}{\tt trifonov@math.sc.edu} \\
\end{center}

\vskip .2 in

\begin{abstract}
We prove that if the least modulus of a distinct covering system is $4$, its largest modulus is at least $60$; also if the least modulus is $3$, the least common multiple of the moduli is at least $120$; finally, if the least modulus is $4$, the least common multiple of the moduli is at least $360$. The constants $60$, $120$, and $360$ are best possible, they cannot be replaced by larger constants.
\end{abstract}
 
\section{Introduction}
A covering system ${\mathcal C}$  is a set of congruences $x \equiv r_i$ (mod $n_i)$, $i = 1,\ldots ,k$,  such that     every integer satisfies at least one of the congruences. Without loss of generality we can assume that ${\displaystyle 1 \leq n_1 \leq \cdots \leq n_k}$. A covering system is distinct if further $1< n_1 < n_2 < \ldots < n_k$. Note that we allow $1$ to be a modulus for a covering system in this paper but do not allow $1$ to be a modulus for a distinct covering system. Throughout the paper we will denote the least modulus $n_1$ of the covering by $m$, the largest modulus $n_k$ by $M$, and the least common multiple of all moduli by $L({\mathcal  C}) = L$. 
For example, 
\begin{equation}\label{eq:m2}
 x \equiv 1 \ (\mbox{mod} \ 2), \  x \equiv 2 \ (\mbox{mod} \ 4), \
 x \equiv 0 \ (\mbox{mod} \ 3), \ x \equiv 4 \ (\mbox{mod} \ 6),  \ \mbox{ and } 
 x \equiv 8 \ (\mbox{mod} \ 12)
\end{equation}
is a distinct covering system with $m=2$, $M=12$, and $L=12$.

Erd\H{o}s \cite{Erd} introduced the use of covering systems in number theory in the 1950s. He constructed a distinct covering system with least modulus $3$ and largest modulus $120$. Erd\H{o}s \cite{Erd} wrote, ``It seems likely that for every $c$ there exists such a system all the moduli
of which are $ > c$."  Proving or disproving this statement became {\it the minimum modulus problem}. For decades many mathematicians believed that indeed, it is possible to construct covering systems with arbitrarily large least modulus.

Swift \cite{Swift} (1954) found a distinct covering with  $m=4$ and later on with $m=6$. This was improved throughout the years by
Churchhouse \cite{ChurchHouse}  with $m=9$ (1968), Krukenberg \cite{Krukenberg}  with $m=18$ (1971), Choi \cite{Choi} with $m =20$ (1971), and Morikawa \cite{Morikawa} with $m=24$ (1981). 
Twenty five years later, Gibson \cite{Gibson} constructed a distinct covering with $m=25$. In 2009, Nielsen \cite{Nielsen} introduced the use of recursion in covering systems and constructed a distinct covering whose smallest modulus is $40$. In the same paper Nielsen wrote, 
``The method further demonstrates some of the difficulty in answering Erd\H{o}s’ minimum modulus
problem, and leads the author to believe that it has a negative solution." Owens \cite{Owens} refined Nielsen's approach and constructed a distinct covering system with minimum modulus $42$.

In 1980, Erd\H{o}s and Graham \cite{ErdGra} investigated systems of congruences with all moduli in $[n, cn]$, where $c > 0$ is a fixed constant. They conjectured that for each $c > 0$ there exists $n(c)$ and $\epsilon(c) > 0$, such that for  each set of congruences with moduli 
$n_1 < \cdots < n_k$ all in $[n , cn]$, the density of the uncovered set is at least $\epsilon (c)$ provided $n$ is sufficiently large, $n \geq n(c)$. 

Erd\H{o}s and Graham's conjecture was proved in 2007 by Filaseta, Ford, Konyagin, Pomerance, and Yu \cite{Fil}. Building on the work of Filaseta et. al., Hough \cite{Hough} made a real breakthrough and solved {\it the minimum modulus problem}. He showed that the minimum modulus in any distinct covering system does not exceed $10^{16}$. 

 Erd\H{o}s and Selfridge posed another famous problem, {\it the odd covering problem}. The problem is to determine whether there exists a distinct covering with all moduli odd integers. Erd\H{o}s was convinced \cite{Erdodd} that such coverings exist and offered  \$25 for a proof that no such covering exists. Selfridge
(recounted from \cite{Fil2}) was convinced that no such covering exists and offered \$2000 for the first example of an odd covering. 

Work of Balister, Bolobas, Morris, Sahasrabudhe, and Tiba brings us the closest to solving the odd covering problem. Balister et al. \cite{Bal0} show that if ${\mathcal C}$ is a distinct covering,  then either $2|L({\mathcal C})$, or $9|L({\mathcal C})$, or $15|L({\mathcal C})$. 
The authors also show that  the least modulus of a distinct covering system does not exceed $606,000$, and that \cite{Bal1} there is no distinct covering system in which all moduli  are odd, squarefree integers.

So, in the last fifteen years several remarkable papers concerning coverings with large minimum modulus appeared. In 1971, Krukenberg \cite{Krukenberg} wrote a Ph.D. thesis where he did an extensive study of covering systems with relatively small minimum modulus and obtained a number of interesting results. 
Unfortunately, none of these results were published in mathematical journals. We outline the main results of Krukenberg's dissertation.  

Krukenberg investigated the following problem. Suppose the least modulus $m$ of a distinct covering system is fixed. What is the least possible value of the largest modulus $M$ of the covering system? The covering system \eqref{eq:m2}  with $m=2$ and $M=12$ has been well-known for many years and Krukenberg constructed a distinct covering system with $m=3$ and $M =36$ and proved:
\begin{theorem}[Krukenberg]\label{Krukenberg1}
(i) If the minimum modulus of a distinct covering system is $2$, then its largest modulus is at least $12$;

\noindent
(ii) If the minimum modulus of a distinct covering system is $3$, then its largest modulus is at least $36$.
\end{theorem}

Krukenberg also found a distinct covering system with $m=4$ and $M=60$. Krukenberg notes that the value of $M=60$ is least possible when $m=4$ and writes ``but this result will not be proved here."  When $m=5$,  Krukenberg constructed a distinct covering system with $M=108$ and conjectured that $108$ is the least possible value of $M$ in this case.

Krukenberg also provided a complete description of all distinct covering systems with least common multiple of the moduli of the form $L = 2^a 3^b$ with $a$ and $b$ positive integers. 

\begin{theorem}[Krukenberg]\label{Krukenberg2}
 Let ${\mathcal C}$ be a distinct covering system with least common multiple of the moduli of the form $L = 2^a 3^b$ with $a$ and $b$ positive integers and least modulus $m$. Then

\noindent
(i) $m \leq 4$;

\noindent
(ii) if $m = 3$,  then $a \geq 3$, and $b \geq 2$;

\noindent
 (iii) if $m=3$  and $a=3$, then $b \geq 3$;

  \noindent
 (iv) there exist coverings with $m = 3$ for each $L\in \{2^4 3^2, 2^3 3^3\}$;
 
 \noindent
 (v) if $m = 4$, then $a \geq 5$ and $b \geq 3$;

\noindent
 (vi) there exist coverings with $m = 4$ for each $L\in\{2^7 3^3, 2^6 3^4, 2^5 3^5\}$ ;

\noindent
 (vii) there is no covering with $m = 4$ and $L \in \{2^6 3^3, 2^5 3^4\}$.
\end{theorem}

Krukenberg also constructed a distinct covering system where all moduli are squarefree integers and the system does not use the modulus $3$. The problem whether there exists a distinct covering system with all moduli squarefree integers and least modulus $3$, is still open.

Finally, for $m$ between $3$ and $18$ with the exception of $m=7$, Krukenberg constructed distinct covering systems with least modulus $m$ while trying to keep the least common multiple $L$ of all moduli  small. Having $L$ small is an advantage. It is much easier to understand the structure of the covering system when $L$ is small and also to modify the covering system to obtain a covering system with different properties. Below is a table comparing $L$ in the systems constructed by Churchhouse \cite{ChurchHouse} to the systems constructed by Krukenberg when $m$ is between $3$ and $9$. 

\vspace{0.15 in}
\noindent
\begin{center} 
$\begin{array}{|c | c| c|}
\hline
m & L {\mbox{ (Churchhouse} )} &  L {\mbox{ (Krukenberg} )}\\ \hline
3 & 2^3\times 3\times 5 & 2^3 \times 3  \times 5  \\ \hline
4 & 2^4 \times 3^2 \times 5 &  2^3 \times 3^2 \times 5\\ \hline
5 & 2^3 \times 3^2 \times 5 \times 7  &  2^5 \times 3^2 \times 5\\ \hline
6 & 2^5 \times 3^2 \times 5 \times 7 &  2^4 \times 3^2 \times 5 \times 7\\ \hline
7 & 2^5 \times 3^3 \times 5 \times 7 &   \\ \hline
8 & 2^4 \times 3^3 \times 5^2 \times 7 &  2^5 \times 3^2 \times 5^2 \times 7\\ \hline
9 & 2^7 \times 3^3 \times 5^2 \times 7 &  2^5 \times 3^3 \times 5^2 \times 7\\ \hline 
\end{array}\ $
\end{center} 

\vspace{0.15 in}
In the table above there is no entry in the third column for $m=7$ since Krukenberg modified the covering with $m=6$ to jump straight to one with $m=8$. 

As a result of the work of Krukenberg, we have an almost complete understanding of distinct covering systems when $m=3$ and $m=4$. In this paper we tie a few loose ends left when $m=3,4$ and lay the groundwork to extend Krukenberg's work to larger $m$. 

Filaseta, Yu, and the second author \cite{Fil3} showed that for each integer $n \geq 3$ there is no distinct covering system with all moduli in the interval $[n, 6n]$. We prove the result with a larger constant.

\begin{theorem}\label{8}
For each integer $n \geq 3$, there is no distinct covering system with all moduli in the interval $[n, 8n]$. 
\end{theorem}

In the fifty years since the Ph.D. thesis of Krukenberg, no proof of Krukenberg's claim that if $m=4$, then $M \geq 60$ has appeared. We supply a proof.

\begin{theorem}\label{m4}
If the minimum modulus in a distinct covering system is $4$, then its largest modulus is at least $60$. 
\end{theorem}

Recall that Churchhouse found a covering with $m=3$ and $L=120$ and Krukenberg found one with $m=4$ and $L=360$. Can one replace the constants $120$ and $360$ by smaller constants? We show that this is not the case.

\begin{theorem}\label{lcm}
(i) If the minimum modulus in a distinct covering is $3$, then the least common multiple of all the moduli is at least $120$;

\noindent
(ii) If the minimum modulus in a distinct covering is $4$, then the least common multiple of all the moduli is at least $360$.
\end{theorem}

The paper is organized as follows. In Section \ref{sec:2}, we introduce new notation which makes analyzing coverings easier, refine an approach of Krukenburg on reducing the  number of congruences in a covering, and prove 
Theorem \ref{8}.
In Section \ref{sec:3}, we introduce another tool, `reduction of a covering' and prove Theorem \ref{m4}. In Section  \ref{sec:4}, we prove Theorem \ref{lcm}. Finally, in Section \ref{open pr}, we formulate some open problems and indicate possible extensions of Krukenberg's work.  

\section{Reducing the number of congruences in a covering system} \label{sec:2}

First, we introduce a notation for congruences which is  convenient when dealing with covering systems. 

Assume we are considering a covering with least common multiple of the moduli 
$L = p_1^{b_1} \cdots p_k^{b_k}$ (unless specified otherwise, $p_k$ will be the $k$th prime number).
Consider the congruence $x \equiv r \ (\mbox{mod} \ n)$, where $n > 1$ has prime factorization $n=p_1^{a_1} \cdots p_k^{a_k}$. For the moment, we suppose $a_l  \geq 1$ for $l=1,\dots,k$.  Next, we find the remainders $r_1,r_2,\ldots,r_k$ when $r$ is divided by $p_1^{a_1},\ldots, p_k^{a_k}$ respectively. 
Let $d_1$ be the base $p_1$ - representation of $r_1$ with its base $p_1$ digits written in {\bf reverse} order. Define similarly, $d_2, \ldots, d_k$. Then,  $x \equiv r \ (\mbox{mod} \ n)$ is written $(d_1 |\ \ldots |\ d_k )$ in our notation.

For example, consider the congruence $ x \equiv 6 \ (\mbox{mod} \ 120)$.  It is equivalent to the system of congruences  $x \equiv 6 \ (\mbox{mod} \ 8)$, $x \equiv 0 \ (\mbox{mod} \ 3)$,
and $x \equiv 1 \ (\mbox{mod} \ 5)$. Thus, for  $ x \equiv 6 \ (\mbox{mod} \ 120)$ we have $120 = 2^3 \times 3 \times 5$, $r_1 = 6$, $r_2 = 0$, $r_3 = 1$, and $d_1 = 011$, $d_2 = 0$, $d_3 = 1$ (since $6$ is $110$ in base $2$). So $x \equiv 6 \ (\mbox{mod} \ 120)$ is written $(011|\ 0|\ 1)$.

A technical note on the above notation is that if we consider a congruence modulo $p_1^{a_1} \cdots p_k^{a_k}$, we make sure that in the new notation we have $a_1$ base $p_1$ digits in the first component, $a_2$ base $p_2$ digits in the second component, and so on. For example, $ x \equiv 0 \ (\mbox{mod} \ 360)$ will be $(000|\ 00|\ 0)$, and not  $(0|\ 0|\ 0)$ (the congruence $(0|\ 0|\ 0)$ is $ x \equiv 0 \ (\mbox{mod} \ 30)$).

The reason we reverse the order of the digits is as follows. Imagine all nonnegative integers organized as a tree with all integers at a vertex at the top, branching to two vertices, one with even integers to the left labeled ($0$) in our notation, and one with odd integers to the right right labeled  ($1$).  Next,  each  of these two vertices branches into two vertices, so we 
get  vertices $(00)$ and $(01)$ on the left, and vertices $(10)$ and $(11)$ on the right. Having the base $2$ digits in reverse order makes it faster to find our path in this tree. 

Furthermore, if  one  or more of the exponents $a_l$ in the factorizatizion $n=p_1^{a_1} \cdots p_k^{a_k}$ is zero, then we put $*$ in the $l$th position of the notation for the congruence. For example,  $ x \equiv 1 \ (\mbox{mod} \ 10)$ is written 
$(1|\  *|\ 1)$. 

Sometimes, it will be possible to write several residue classes in a more compact way. For example, suppose that at a certain stage of constructing a covering, the uncovered set consists of the residue classes $x \equiv 0 \ (\mbox{mod} \ 6)$ and 
$x \equiv 4 \ (\mbox{mod} \ 6)$. In this case, we will denote the uncovered set by $(0|\  0,1)$. 

Finally, for brevity we truncate trailing $*$s. For example, if $L = 60$, the congruence $ x \equiv 0 \ (\mbox{mod} \ 2)$ will be written as $(0)$ rather than $(0|\ *|\ *)$. 

For a final example on this notation, let us analyze the distinct covering system given in \eqref{eq:m2}. 
The first two congruences are $(1)$ and $(01)$ leaving a congruence class modulo $4$, namely  $(00)$, uncovered. We split it into three classes modulo $12$, namely 
$(00|\ 0,1,2)$ which is our way of writing the three congruences given by $(00|\ 0)$,  $(00|\ 1)$,  and $(00|\ 2)$. 
We cover $(00|\ 0)$ by a congruence modulo $3$, $(*|\ 0)$; we cover $(00|\ 1)$ by a congruence modulo $6$, $(0|\  1)$; finally, we cover $(00|\ 2)$ by a congruence modulo $12$, $(00|\ 2)$. 

We refer to the representation of a residue class we just introduced as a coordinate representation.
This  notation  is in line with the geometric approach to covering systems of Simpson and Zeilberger \cite{SimZeil}. In the case when $L$ is squarefree, congruences correspond to points and hyperplanes in a certain $k$ dimensional box.

If $p$ is a prime, $a$ is a nonegative integer, and $n$ is a positive integer, then $p^a || n$ will mean that $p^a | n$ and $p^{a+1} \nmid  n$. 

Next, we define two operations on  residue classes - splitting modulo $p$ and reducing  modulo $p$. 

Assume that $p$ is a prime, $a$ is a nonegative integer,  $n$ is a positive integer,  and $p^a || n$. Splitting the residue class $r \ (\mbox{mod} \ n) $ modulo $p$ means that we replace it by $p$ residue classes modulo $np$  (`fibers') by consecutively appending the base-$p$ digits  $0,1,\dots , p-1$  in the position corresponding to 
$p^{a+1}$ in  the coordinate representation of the residue class. We denote the $l$th `fiber' described above by $(r \ (\mbox{mod} \ n))_{p,l} $.
For example, if we split $(1|\ 1|\ 4)$ modulo $3$, we obtain the  three `fibers'  $(1|\ 10,11,12 |\ 4)$. 

Similarly, assume that $p$ is a prime, $a$ and $n$ are positive integers, and $p^a || n$. Reducing the residue class $r \ (\mbox{mod} \ n) $ modulo $p$ means that we delete the base-$p$ digit in the position corresponding to 
$p^a$ in the coordinate representation of the residue class. 
For example, if we reduce $(0|\ 21|\ 34)$ modulo $5$ we get $(0|\ 21|\ 3)$ .

Our first tool is the following lemma which builds on ideas of  Krukenberg \cite{Krukenberg}.

\begin{lemma} \label{lem:1}
Let ${\mathcal C}$ be a covering system with least common multiple of the moduli $L$. Assume $p^a || L$ for some prime $p$ and a positive integer $a$. Denote by ${\mathcal C_0}$ the subset of congruences in ${\mathcal C}$ whose moduli are not divisible by $p^a$; also, let 
${\mathcal C_1}$ be the subset of congruences in ${\mathcal C}$ whose moduli are  divisible by $p^a$.

Next, for $l=0,\ldots,p-1$, define $B_l$ as the subset of congrences in ${\mathcal C_1}$ whose modulus has base-$p$ digit corresponding to $p^a$ (in coordinate notation) equal to $l$. 

Finally, let $D_l$ be the set of congruences in $B_l$ reduced modulo $p$. 

Then, one can replace the congruences in ${\mathcal C_1}$ by ${\displaystyle D = \bigcap_{l=0}^{p-1} D_l}$ and we will still have a covering;  that is, ${\mathcal C_0} \cup D$ is a covering system. 

\end{lemma}

To clarify, what we do is sort the congruences with moduli divisible by $p^a$ by the base-$p$ digit corresponding to $p^a$ in `bins' $B_l$. Next, we delete the base-$p$ digits corresponding to $p^a$ from all congruences in the bins. 
Finally, we take the intersection of the union of the reduced congruences in each bin.  Note  that the intersection of  unions of sets can be written as a union of intersections. Also,  the intersection of the sets covered by several congruences is either an empty set   or the set covered by a single congruence with modulus the least common multiple of the moduli of the congruences we intersect.  
The claim is that we can replace the congruences in ${\mathcal C_1}$ by the congruences  we obtain by the process described above.

For example if  we apply Lemma \ref{lem:1} with $p=3$ to the covering in \eqref{eq:m2}, the `bins' are $B_0 = \{ (*|\ 0) \}$, $B_1 = \{ (0|\ 1) \}$, and $B_2 = \{ (00|\ 2) \}$. Reducing modulo $3$ we get
$D_0 = \{ (*|\ *) \}$, $D_1 = \{ (0|\ *) \}$, and $D_2 = \{ (00|\ *) \}$. So,  $D  = \{ (00|\ *) \}$. We claim that  replacing the congruences with moduli $3$, $6$, and $12$ by a single congruence modulo $4$ still leaves us with a covering. 
Indeed, $(1)$, $(01)$, and $(00)$ is a covering.

\begin{proof}
Let $R$ be the set uncovered by the conguences in ${\mathcal C_0}$. Note that the least common multiple of the moduli in ${\mathcal C_0}$ divides $L_1=L/p$. Therefore, $R$ can be expressed as a union of residue classes modulo $L_1$ and 
$p^{a-1} || L_1$.  

Let $r \ (\mbox{mod} \ L_1) $ be one of the uncovered residue classes in $R$.  
 We split it modulo $p$. Consider the `fiber'  $(r \ (\mbox{mod} \ L_1))_{p,0}$. It does not satisfy any of the congruences in ${\mathcal C_0}$ or in `bins' $B_1, \ldots, B_{p-1}$. We say that a set of congruences $C$ covers a certain set of integers $S$ if every integer in $S$ satisfies at least one congruence in $C$. Then, since ${\mathcal C}$ is a covering, 
the congruences in `bin' $B_0$ cover $(r \ (\mbox{mod} \ L_1))_{p,0}$. 

Reducing modulo $p$ we get that the congruences in `bin' $D_0$ cover $r \ (\mbox{mod} \ L_1) $. Since this is true for each residue class in $R$, we get that the congruences in `bin' $D_0$ cover $R$. 

Similarly, we get that the congruences in `bin' $D_l$ cover $R$ for each $l=1,\dots,p-1$. Therefore, ${\displaystyle D =  \bigcap_{l=0}^{p-1} D_l}$ covers $R$, so we can replace  ${\mathcal C_1}$ by $D$ and will still have a covering.
\end{proof}

\begin{corollary} \label{cor:1}
Let ${\mathcal C}$  be a covering such that $p^a | L$ for some prime $p$ and integer $a \geq 1$. Suppose that there are $k$ congruences in ${\mathcal C}$ whose moduli are divisible by $p^a$. Then,  if $k < p$, we can discard from ${\mathcal C}$ all congruences 
whose moduli are divisible by $p^a$, and will still have a covering.
\end{corollary}

\begin{proof}
First, we justify that we need only consider the case $a=\nu_p(L)$, where $\nu_p(m)$ for $m \in \mathbb{N}$ is the integer for which $p^{\nu_p(m)}|| m$. Suppose we have established the result for the case $a=\nu_p(L)$. With $k$ as stated in the corollary, $k$ is an upper bound on the number of congruences in ${\mathcal C}$ with moduli divisible by $p^j$ for all $j \in \mathbb{Z}$ with $a \leq j \leq \nu_p(L)$. Then by applying the corollary for $a$ replaced by $\nu_p(L)$ to obtain a new covering and then applying the corollary over and over again, one arrives at the covering ${\mathcal C}$ with all congruences having moduli divisible by $p^a$ removed, proving the corollary. So we suppose now $a=\nu_p(L)$.

Now, let $a=\nu_p(L)$. Since $k<p$, we see that there is an $l \in \{0, 1, \ldots, p-1\}$ in Lemma~\ref{lem:1} such that $B_l \neq \emptyset$. Therefore, $D_l$ and, hence, $D$ is $\emptyset$. The corollary now follows from Lemma~\ref{lem:1}.
\end{proof}

\begin{corollary} [Krukenberg]  \label{cor:2}
Let ${\mathcal C}$  be a distinct covering  with all moduli in the interval $[c,d]$.  If $p$ is a prime and $a$ is a positive integer such that $p^a (p+1) >  d$, then we can discard all congruences whose moduli are multiples of $p^a$ and still have a covering.
\end{corollary}

\begin{proof}
First, since $p^a(p+1) < 2p^{a+1}$, there is at most one multiple of $p^{a+1}$ in  $[c,d]$ so by Corollary \ref{cor:1} we can discard the congruence with modulus $p^{a+1}$ (if there is one). This leaves us with at most $p-1$ multiples of $p^a$ in $[c,d]$, namely 
$p^a \cdot 1, \ldots , p^a \cdot (p-1)$. By applying Corollary \ref{cor:1} again, we can discard all moduli divisible by $p^a$. 
\end{proof}

\begin{corollary} [Krukenberg]  \label{cor:3}
Let ${\mathcal C}$  be a covering such that $p^a || L$ for some prime $p$ and integer $a \geq 1$. Let ${\mathcal C_1}$  be the subset of ${\mathcal C}$ consisting of congruences  whose moduli are divisible by $p^a$. Suppose $|{\mathcal C_1 }| = p$ and the moduli of the congruences in 
${\mathcal C_1}$ are $p^a m_1, \ldots, p^a m_p$. Then, 

\noindent
(i) one can replace the congruences in ${\mathcal C_1}$ by a single congruence with modulus
$${\displaystyle p^{a-1}  \lcm (m_1,\ldots, m_p)}$$ and the resulting set will still be a covering.

\vskip 5pt
\noindent
(ii) if two of the above $p$ congruences are in the same class modulo $p^a$ we can discard all $p$ congruences and the resulting set will still be a covering.

\end{corollary}

\begin{proof}
(i) Again we use  Lemma \ref{lem:1}. Now,  we split the $p$ congruences into $p$ `bins' $B_l$, with $l \in \{0, 1, \ldots, p-1\}$,  as in Lemma \ref{lem:1}. The only case when $D \neq \emptyset$ is when there is exactly one congruence in each `bin' and when the system of the $p$ congruences reduced modulo $p$ with moduli $p^{a-1} m_1, \ldots, p^{a-1} m_p$ has a solution. 
By a generalization of the Chinese remainder theorem, if a finite system of congruences has a solution, the system of congruences is equivalent to a single congruence whose modulus is the least common multiple of the congruences in the finite system. Lemma \ref{lem:1} now implies \emph{(i)}.

\vskip 5pt
\noindent
(ii) Here we note that since a `bin' contains two or more congruences, at least one of the remaining `bins' will be empty, so $D = \emptyset $ in this case. Then the conclusion of Lemma \ref{lem:1} implies \emph{(ii)}.
\end{proof}

Next, we define a {\it   minimal covering system}. A minimal covering system ${\mathcal C}$ is a covering such that no proper subset of ${\mathcal C}$ is a covering system. Clearly, by discarding one by one redundant congruences, after a finite number of steps, any finite covering system can be reduced to a minimal covering system in at least one way. 

Next, we use Lemma \ref{lem:1} to prove Theorem \ref{8} under the assumption that Theorem \ref{m4} holds.

\begin{proof}[Proof of Theorem \ref{8} assuming Theorem \ref{m4}.]
Assume that for some integer $m \geq 3$ there is a distinct covering ${\mathcal C}$ with all moduli in the interval $[m,8m]$. Let ${\mathcal C_m}$ be a minimal covering which is a subset of ${\mathcal C}$. Consider the least common multiple $L$ of the moduli of the congruences in 
${\mathcal C_m}$. By Corollary \ref{cor:1}, if $p^a | L$ for some prime $p$ and a positive integer $a$, then the interval $[m,8m]$ contains at least $p$ multiples of $p^a$ that are not multiples of $p^{a+1}$. Since one of every $p$ consecutive multiples of $p^a$ is divisible by $p^{a+1}$, we deduce that the interval $[m,8m]$ contains at least $p+1$ multiples of $p^a$.

Denote by ${\mathcal M}\subseteq [m, 8m]$ the set of moduli from the congruences in ${\mathcal C_m}$.
Let $p \geq \sqrt{7m+1}$  be a prime. The number of multiples of  $p$ in the interval $[m, 8m]$ is 
$$n_p :=\left  \lfloor {\frac{8m}{p}}  \right \rfloor - \left  \lfloor {\frac{m-1}{p}} \right \rfloor = {\frac{7m+1}{p}} - \left \{  {\frac{8m}{p}} \right  \}  +
 \left \{  {\frac{m-1}{p}} \right  \},$$ where $\{ x \}$ denotes the fractional part of $x$. 
Since, for each $x$, $0 \leq \{ x \} < 1$, we get $$n_p <{\displaystyle  \frac{7m+1}{p} }  +1 \leq \sqrt{7m+1} + 1 \leq  p+1.$$ Thus, for each $p \geq \sqrt{7m+1}$,  there are less than $p+1$ multiples of $p$ in the interval $[m, 8m]$. Therefore, if $n$ is a modulus of one of the congruences in  ${\mathcal C_m}$ (that is $n \in {\mathcal M}$), then all the prime divisors of $n$ are 
less than $\sqrt{7m+1}$. Since, the density of integers covered by a  congruence modulo $n$ is $1/n$  and ${\mathcal C_m}$ is a covering, we get 

\begin{equation} \label{m3}
  \sum_{\substack{ m \leq n \leq 8m, \\ P(n) < \sqrt{7m+1}}} \frac{1}{n} \geq \sum_{n \in {\mathcal M}}  \frac{1 }{n} \geq 1,
\end{equation}
where $P(n)$ denotes the largest prime divisor of $n$.

Let $$S_m = \sum_{n \in {\mathcal M}} \frac{1}{n} \quad  \text{ and } \quad T_m = \sum_{ \substack{m \leq n \leq 8m, \\ P(n) < \sqrt{7m+1}}} \frac{1}{n}.$$ 

We checked by direct computation and by using the inequality $T_{m-1} \leq T_m + \frac{1}{m-1}$ that $T_m < 1$ for all $m \in [51,606000]$. Details on this computation are in  Appendix 1 to the paper. Since Ballister et al. \cite{Bal0} showed that the minimum modulus of a distinct covering system does not exceed 
$606000$,   Theorem \ref{8} holds when $m \geq 51$.  

Also, since Krukenberg showed that there is no distinct covering system with moduli in $[3,35]$,  Theorem \ref{8} holds when $m=3$, and $4$. 

Furthermore, given Theorem \ref{m4}, there is no distinct covering system with moduli in $[4,59]$; therefore Theorem \ref{8} holds when $m=5, 6,$ and $7$. 

There are twenty occasions when $m \in [8,50]$ and $T_m \geq 1$. They are shown below.

\begin{center}
\vspace{0.15 in} \noindent
$\begin{array}{|c|c|c|c|c|c|c|}
\hline
m & 8 & 9 & 10 & 11 & 12 \\ \hline
T_m & 1.265378\ldots & 1.168553\ldots & 1.083275\ldots & 1.007526\ldots & 1.029053\ldots \\ \hline \hline
m & 15 & 16 & 17 & 18 & 20 \\ \hline 
T_m & 1.105901\ldots & 1.071248\ldots & 1.023731\ldots & 1.037819\ldots & 1.008476\ldots \\ \hline \hline
m & 21 & 22 & 23 & 24 & 25 \\ \hline 
T_m & 1.138447\ldots & 1.108141\ldots & 1.073737\ldots & 1.084236\ldots & 1.062850\ldots \\ \hline \hline
m & 26 & 47 & 48 & 49 & 50 \\ \hline 
T_m & 1.027658\ldots & 1.028747\ldots & 1.036628\ldots & 1.023508\ldots & 1.01063\ldots \\ \hline
\end{array}$
\end{center}

\vspace{0.15 in} 
So far, we have used Corollary \ref{cor:1} only with $a=1$. Next we use Corollary \ref{cor:2} for all $a \geq 1$.

Define $$L_m=
\begin{cases}
720  = 2^4\cdot 3^2 \cdot 5 \text{ if } m=8 \\
5040  = 2^4\cdot 3^2 \cdot 5 \cdot 7  \text{ if } m\in \{9, 10, 11\} \\
10080  = 2^5\cdot 3^2 \cdot 5 \cdot 7 \text{ if } m=12 \\
30240  = 2^5\cdot 3^3 \cdot 5 \cdot 7 \text{ if } m\in \{15, 16\} \\
332640  = 2^5\cdot 3^3 \cdot 5 \cdot 7 \cdot 11 \text{ if } m\in\{17, 18\} \\
1663200 = 2^5\cdot 3^3 \cdot 5^2 \cdot 7 \cdot 11 \text{ if } m\in \{20, 21, 22\} \\
21621600 = 2^5\cdot 3^3 \cdot 5^2 \cdot 7 \cdot 11 \cdot 13 \text{ if } m=23 \\
43243200 =2^6\cdot 3^3 \cdot 5^2 \cdot 7 \cdot 11 \cdot 13 \text{ if } m\in \{24, 25, 26\} \\
2205403200 =2^6\cdot 3^4 \cdot 5^2 \cdot 7 \cdot 11 \cdot 13 \cdot 17\text{ if } m=47 \\
83805321600 =2^7\cdot 3^4 \cdot 5^2 \cdot 7^2 \cdot 11 \cdot 13 \cdot 17 \cdot 19 \text{ if } m=48 \\
586637251200 = 2^7\cdot 3^4 \cdot 5^2 \cdot 7^2 \cdot 11 \cdot 13 \cdot 17 \cdot 19 \text{ if } m\in\{49, 50\}. 
\end{cases}$$
Using Corollary \ref{cor:2}, one checks directly that $L$ divides $L_m$ for each $$m\in\{8, 9, 10, 11, 12, 15, 16, 17, 18, 20, 21, 22, 23, 24, 25, 26, 47, 48, 49, 50\}.$$ Then for such $m$ we have 
\begin{equation} \label{eq3}
\sum_{n\in M} \frac{1}{n} \leq \sum_{\substack{d \mid L_m \\ m \leq d \leq 8m}} \frac{1}{d} < 1,
\end{equation}
where the last inequality is done by a direct computation. As \eqref{eq3} contradicts the second inequality in \eqref{m3}, the proof is complete.
\end{proof}

A natural question is whether one can replace the constant $8$ in Theorem \ref{8} by a larger constant? If we try to prove 
Theorem \ref{8} with a constant $9$, the estimate of $T_m$ for large $m$ is similar to what we did above. However, there will be   $97$  
values of $m$ for which $T_m > 1$ and dealing with these would make the proof of the theorem much longer.  Moreover, dealing with some of the smaller exceptions $m$ 
would require more intricate approach than the one in the proof of Theorem~\ref{8}.

\section{Reduction of a covering}\label{sec:3}

Our second tool is {\it reduction of a covering}. We start with an example. 

Consider the covering \eqref{eq:m2}. Let $a \in \{ 0,1,2 \}$. Since  \eqref{eq:m2} is a covering, the residue class $3m + a$ is covered by the congruences in \eqref{eq:m2}.

Substituting $3m + a$ for $x$ in each of the congruences of \eqref{eq:m2}  and solving for $m$ we get a new covering system. When $a=0$ we get 
$m \equiv 1 \ (\mbox{mod} \ 2)$,  $ m \equiv 2 \ (\mbox{mod} \ 4)$, 
$ m \equiv 0 \ (\mbox{mod} \ 1)$,  and two congruences have no solution; when $a=1$ we obtain $m \equiv 0 \ (\mbox{mod} \ 2)$,  $ m \equiv 3 \ (\mbox{mod} \ 4)$, 
$ m \equiv  1 \ (\mbox{mod} \ 2)$,  and two congruences have no solution; finally, when $a=2$ we have\\ $m \equiv 1 \ (\mbox{mod} \ 2)$,  $ m \equiv 0 \ (\mbox{mod} \ 4)$, 
$ m \equiv 2 \ (\mbox{mod} \ 4)$,  and two congruences have no solution.

In general, let ${\mathcal C}$ be a covering system  and let $p$ be a prime. Let ${\mathcal C_0}$ be the subset of ${\mathcal C}$ of congruences whose moduli are not divible by $p$. Let ${\mathcal M_0}$ be the list of the moduli of the congruences in 
${\mathcal C_0}$. Similarly, let ${\mathcal C_1}$ be the subset of ${\mathcal C}$ of congruences whose moduli are  divisible by $p$, and  let ${\mathcal M_1}$ be the list of the moduli of the congruences in 
${\mathcal C_1}$. 

To reduce the covering modulo $p$ for each $a \in \{ 0,1,\ldots, p-1 \}$ we substitute $pm+a$ for $x$ in each of the congruences of
${\mathcal C}$ and solve for $m$ to get a new covering. This way we end up with $p$ coverings in which each modulus in  ${\mathcal M_0}$ is  used in all $p$ coverings. However, if $m$ is a modulus in ${\mathcal M_1}$ it gets replaced by 
$m/p$ and it is used in just one of the $p$ coverings. 

Indeed, if $x \equiv r \ (\mbox{mod} \ n)$ is a congruence in ${\mathcal C_0}$ (so $p \nmid n$), substituting $mp+a$ for $x$ and solving for $m$, we get $m \equiv p^{-1}(r-a) \ (\mbox{mod} \ n)$.

However, if $x \equiv r \ (\mbox{mod} \ n)$ is a congruence in ${\mathcal C_1}$ (so $p  | n$), substituting $mp+a$ for $x$, we get the congruence $mp+a \equiv r \ (\mbox{mod} \ n)$. The last congruence has a solution if and only if 
$r \equiv a \ (\mbox{mod} \ p)$, in which case we get $m \equiv (r-a)/p \ (\mbox{mod} \ n/p)$.

Next, we say that two congruences are in the same class modulo a positive integer $q$, if the integers covered by the congruences all belong to one class modulo $q$. In other words, if the two congruences are 
$x \equiv r_1 \ (\mbox{mod} \ n_1)$ and $x \equiv r_2 \ (\mbox{mod} \ n_2)$, we say that they are in the same class modulo $q$ if $q | n_1$, $q | n_2$, and $r_1  \equiv r_2 \ (\mbox{mod} \ q)$.

Assume two congruences $x \equiv r_1 \ (\mbox{mod} \ n_1)$ and $x \equiv r_2 \ (\mbox{mod} \ n_2)$, both in ${\mathcal C}$ are in the same class modulo $q$ with $p \nmid q,n_1,n_2$. After reduction modulo $p$, 
we get $mp \equiv ( r_1-a) \ (\mbox{mod} \ n_1)$ and  $mp \equiv ( r_2-a) \ (\mbox{mod} \ n_2)$. Since $q | n_1$, $q | n_2$, and $p \nmid q,n_1,n_2$,   the reduced congruences are still in the same class modulo $q$.

To summarize, we showed that the following lemma holds. 

\begin{lemma} \label{lem:2}
Let ${\mathcal C}$ be a covering system  and let $p$ be a prime. Let ${\mathcal C_0}$ be the subset of ${\mathcal C}$ of congruences whose moduli are not divisible by $p$. Let ${\mathcal M_0}$ be the list of the moduli of the congruences in 
${\mathcal C_0}$. Similarly, let ${\mathcal C_1}$ be the subset of ${\mathcal C}$ of congruences whose moduli are  divisible by $p$, and  let ${\mathcal M_1}$ be the list of the moduli of the congruences in 
${\mathcal C_1}$. 

Reducing the covering ${\mathcal C}$ modulo $p$ produces $p$ coverings where 

\noindent
(i) each modulus in  ${\mathcal M_0}$ is  used in each of the $p$ coverings but each modulus $n$ in ${\mathcal M_1}$ is replaced by $n/p$ and is used in just one of the $p$ coverings, 

and

\noindent
(ii) if   two congruences in ${\mathcal C_0}$ are in the same class  modulo a positive integer $q$, where $p \nmid q$,  then after reduction they are in the same class modulo $q$ in each of the  $p$ coverings.

\end{lemma}

Next, we prove that there is no distinct covering system with moduli in the interval $[4,59]$ using a proof by contradiction. Our proof proceeds as follows. First, we will use Corollary \ref{cor:2} and Corollary \ref{cor:3} to reduce the list of possible moduli in the covering to $17$ integers. Next, we reduce the covering modulo $3$. We explore all ways in which it is possible to construct the three coverings from Lemma \ref{lem:2} which satisfy condition (i) of the lemma. 
It turns out, this can be done in two ways. In both cases, we obtain a contradiction by showing condition (ii)  of Lemma \ref{lem:2} with $q=5$ is violated. 

The reason the details of the proof of Theorem \ref{m4} are somewhat complicated is that using moduli in $[4,59]$ one can get very close to a covering; more precisely, one  can   cover $179$ out of $180$ classes modulo $180$. In Appendix 2 we give an example of a distinct 
covering system using congruences with moduli in $[4,56]$ and a congruence modulo $180$. 

\begin{proof}
Assume that there exists a distinct covering   ${\mathcal C}$ with all moduli in $[4,59]$. Since every covering contains a subset which is a minimal covering, without loss of generality we can assume that 
${\mathcal C}$ is a minimal covering. 
Let ${\mathcal M}$ be the set of the moduli of the congruences in ${\mathcal C}$. Also, let   $L$ be the least common multiple of the moduli   in  ${\mathcal M}$. 

By Corollary \ref{cor:2}, if $p^a(p+1) > 59$ for some prime $p$ and a positive integer $a$, then $p^a$ does not divide any modulus in ${\mathcal M}$, so $p^a \nmid L$.
Since, $2^5 \cdot 3 > 59$, $3^3 \cdot 4 > 59$, $7^2 \cdot 8 > 5^2 \cdot 6 > 59$, and $p(p+1) > 59$ for $p \geq 11$, we get $L | \left ( 2^4 \cdot 3^2 \cdot 5 \cdot 7 \right )$. 
Therefore,  $${\mathcal M} \subseteq \{ 4,8,16,  6,12,24,48,  9,18,36, 5,10,20,40, 15,30,45, 7,14,21,28,35,42,56 \}.$$
Without loss of generality we can assume that $${\mathcal M} = \{ 4,8,16,  6,12,24,48,  9,18,36, 5,10,20,40, 15,30,45, 7,14,21,28,35,42,56 \}.$$
Indeed, for each modulus $m$ is  in the displayed set above which is not in ${\mathcal M}$, we simply add a congruence  $x \equiv 0 \ (\mbox{mod} \ m)$ to ${\mathcal C}$.

When analyzing or constructing a covering using a given set of moduli, following Krukenberg \cite{Krukenberg}, Nielsen \cite{Nielsen}, Balister et. al. \cite{Bal1}, we use the moduli in increasing order of arithmetic complexity. For example, 
above we first list powers of $2$, next, powers of $2$ times $3$, etc. and we do not introduce moduli which are multiples of not yet used   prime $p$, until all moduli with prime divisors less than $p$ are used.

Next, by Corollary \ref{cor:3} we can replace the seven congruences in ${\mathcal C}$ with moduli divisible by $7$ by a single congruence modulo $120 =\mbox{ lcm}(1,2,3,4,5,6,8)$ and still have a covering. 
Also, by Corollary \ref{cor:3} we can replace the two congruences with moduli $16,48$ by a single congruence modulo $24$ and still have a covering. Denote the resulting covering by 
${\mathcal C'}$ and denote the {\it  list}  of the moduli of the congruences in ${\mathcal C'}$ by ${\mathcal M'}$. Then, ${\mathcal M'}$ is   the list 
$[ 4,8,  6,12,24, 24,  9,18,36, 5,10,20,40, 15,30,120, 45 ]$, where $[\ldots ]$ is used to emphasize that ${\mathcal M'}$ is a list. Note that $24$ appears twice in the last list meaning that we   have two congruences modulo $24$ in ${\mathcal C'}$.

This concludes the first part of the proof. 

Next, we reduce the covering ${\mathcal C'}$ modulo $3$ and obtain three coverings, say ${\mathcal C_0'}$,  ${\mathcal C_1'}$,  ${\mathcal C_2'}$, whose  moduli are  ${\mathcal M_0'}$,  ${\mathcal M_1'}$,  ${\mathcal M_2'}$, respectively. 
By Lemma \ref{lem:2} the moduli in ${\mathcal M_0} = \{ 4,8,5,10,20,40 \}$  can be used in all three coverings, and each modulus  in\\ ${\mathcal M_1} = [ 2*,4*,8*,8*,3*,6*,12*,5*,10*,40*,15* ]$ can be used in just one of the coverings. 
We will use * behind a modulus to indicate that it can be used in at most   one of  the three coverings  ${\mathcal C_0'}$,  ${\mathcal C_1'}$,  ${\mathcal C_2'}$.

After relabeling, we can assume that $2*$ is in ${\mathcal M_0'}$. 

Moreover, by Corollary \ref{cor:1}, we can take congruences with moduli $3*,6*,12*,15*$ in only one of ${\mathcal C_0'}$,  ${\mathcal C_1'}$,  ${\mathcal C_2'}$. (If we have just one or two congruences with moduli divisible by $3$ in a covering, we can discard them.) 

It is relatively easy to see that we can take $3*,6*,12*,15*$ to be   in ${\mathcal M_1'}$ or ${\mathcal M_2'}$.  Indeed, $2*,4,8$ are already in ${\mathcal M_0'}$.  If $8*$ is also in ${\mathcal M_0'}$, no additional congruences are needed to 
construct the covering  ${\mathcal C_0'}$. Note that  the congruences with moduli $3*,6*,12*$  can cover a congruence class modulo $4$ (and some integers outside of it). If $8*$ is in  ${\mathcal M_1'}$ or  ${\mathcal M_2'}$ and $3*,6*,12*,15*$ are in ${\mathcal M_0'}$, we can swap the congruences with moduli $3*,6*,12*,15*$ and the congruence modulo $8*$ and we
will still have three coverings.

Thus,  after relabeling we can assume that $3*,6*,12*,15*$ are  in  ${\mathcal M_1'}$.

Next, we analyze how to allocate the moduli $5*,10*,40*$.
We claim that ${\mathcal M_2'}$ contains at least one of the moduli $5*,10*,40*$. Otherwise, 
$${\mathcal M_2'} \subseteq [ 4,8,5,10,20,40,  4*,8*,8* ].$$ Using Corollary \ref{cor:1} we discard any of $5,10,20,40$  from ${\mathcal M_2'}$, leaving us with the impossible task of contructing a covering using only congruences with moduli $4,8,4*,8*,8*$. 

 So, we allocated $5$ out of $11$ moduli in ${\mathcal M_1}$ and have partial information about  three of the remaining six moduli. 
It is possible to allocate the   moduli $4*,8*,8*,5*,10*,40*$ and to construct the three coverings. 
However, we will show that it cannot be done without violating condition (ii) of Lemma \ref{lem:2}. 
To this end, we  consider two cases. 
\vskip 5pt
\noindent
{\bf Case I:} The congruences with moduli $20$ and $40$ in ${\mathcal C'}$ are not  in the same class modulo $5$. 
\vskip 5pt
By Lemma \ref{lem:2} the congruences with moduli $20$ and $40$ are not  in the same class modulo $5$ in  ${\mathcal C_0'}$,  ${\mathcal C_1'}$,  and ${\mathcal C_2'}$, as well. 

First, note that currently  ${\mathcal M_0'}$ contains $\{2*, 4,8,5,10,20,40 \}$ which is not sufficient to construct a covering. (Discard $5,10,20,40$ using Corollary \ref{cor:1} and we are left only with moduli $2*,4,8$.) We assign 
$40*$ to ${\mathcal M_0'}$, so that covering is possible with moduli from ${\mathcal M_0'}$. All the remaining five moduli $4*,8*,8*,5*,10*$ are divisors of $40*$, so the remaining two coverings cannot benefit from us assigning a different modulus instead of $40$ to ${\mathcal M_0'}$.

Next, we concentrate on allocating $5*,10*$.  We proved above that  ${\mathcal M_2'}$ contains at least one of the moduli $5*,10*$ (since $40*$ is already allocated to ${\mathcal M_0'}$). 
There are two subcases. 
\vskip 5pt
\noindent
{\bf Subcase A:} Exactly one of the moduli $5*,10*$ is in ${\mathcal M_2'}$.
\vskip 5pt 
In this subcase, ${\mathcal M_2'} \subseteq [ 4,8,5,10,20,40,  4*,8*,8*, 5*]$, or\\ ${\mathcal M_2'} \subseteq [ 4,8,5,10,20,40,  4*,8*,8*, 10*]$.  By Corollary \ref{cor:3} we can replace 
$5,10,20,40,$ and one of $5*$ and $10*$ by a congruence modulo $8$. Now, we need to construct a covering using moduli $4,8,8$ and some of $4*,8*,8*$. The covering  ${\mathcal C_2'}$ can be completed only if all three moduli $4*,8*,8*$ are in 
${\mathcal M_2'}$.  

Without loss of generality we may assume we are left with  $${\mathcal M_1' } = [ 4,8,5,10,20,40, 3*,6*,12*,15*, 5* ].$$ We finish this subcase by proving the following lemma.

\begin{lemma} \label{lem:3}
There is no covering with congruences whose  moduli form the list\\ $ [ 4,8, 3,6,12,5,5, 10,20,40,15 ]$, such that the congruence with modulus $20$ and the congruence with modulus $40$  are not in the same class modulo $5$.
\end{lemma}
 
\begin{proof}
We assume that there is such a covering and use Lemma \ref{lem:2} to reduce the covering modulo $3$. We get three coverings where each has moduli in the list  $ [ 4,8, 5,5, 10,20,40 ]$  and each of the the moduli in 
$[ 1*, 2*, 4*, 5* ]$ is used in exactly one covering. Consider one of the three coverings, say  ${\mathcal C'' }$ which does not use moduli $1*$ and $2*$. The moduli of ${\mathcal C'' }$  are in the list $ [ 4,4*, 8, 5,5, 5*, 10,20,40 ]$. Next we apply Lemma \ref{lem:1} 
with $p=5$ to the congruences with moduli $5,5,5*,10,20,40$. After reduction modulo $5$ we need to place the reduced congruences with moduli $1,1,1,2,4,8$ in five `bins' and take the intersection of the cogruences in the five `bins'. 
Consider the `bin' which does not contain $1,1,1,2$. This `bin' contains the congruence modulo $4$ or the congruence modulo $8$ both not both, since $4$ and $8$ are not in the same `bin' (the congruences with moduli $20$ and $40$ are not in the same class modulo $5$).  So, $D$ is inside a residue  class modulo $4$.  Thus, we can replace the congruences with moduli $5,5,5*,10,20,40$ by a single congruence modulo $4$. This is a contradiction because it requires building a covering with congruences with moduli $4,4,4,8$.
\end{proof}

\noindent
{\bf Subcase B:} Both moduli $5*,10*$ are in ${\mathcal M_2'}$.
\vskip 5pt
Here $$ {\mathcal M_2'}\subseteq [ 4,8,5,10,20,40, 4*,8*,8*, 5*,10*].$$ We apply  Lemma \ref{lem:1} with $p=5$ to the congruences with moduli $5,10,20,40,5*,10*$. Proceeding word by word like  in the proof of Lemma \ref{lem:3} 
we get that we can replace the congruences with moduli $5,10,20,40,5*,10*$ by a single cogruence modulo $4$. Now, we need to construct a covering using moduli $4,8,4$ and some of the moduli $4*,8*,8*$. 
This can be done only if $4*$ and at least one $8*$ are in  ${\mathcal M_2'}$.

This leaves  $${\mathcal M_1' } = [ 4,8,5,10,20,40, 3*,6*,12*,15*, 8* ].$$ Next, by using Corollary \ref{cor:3} we replace the congruences with moduli $5,10,20,40, 15*$ by a single congruence modulo $24$. So, we need to construct a covering using moduli $4,8, 8, 3,6,12,24$. Assume there is such a covering and reduce it modulo $3$. We get three new coverings with common moduli $4,8,8$ and moduli to be used by just one of the three coverings: $1*,2*,4*,8*$. 
Consider the covering which does not contain $1*,2*$. The moduli of the congruences of this covering are at most $4,8,8,4*,8*$, which is a contradiction.  

This completes Case I. 
\vskip 5pt
\noindent
{\bf Case II:} The congruences with moduli $20$ and $40$ in ${\mathcal C'}$ are in the same class modulo $5$. 
\vskip 5pt
First, we claim that ${\mathcal M_2'}$ contains at least two of the moduli $5*,10*,40*$. Assume otherwise. Then, ${\mathcal M_2'}$ contains at most the moduli $4,8,5,10,20,40, 4*,8*, 8*$ and one of the moduli $5*,10*, 40*$.
By Corollary \ref{cor:3} we can discard from ${\mathcal C_2'}$ the congruences with moduli $5,10,20,40$ and the congruence with modulus $5*$, $10*$, or $40*$ since since of two these congruences are in the same class modulo $5$. This leaves us at most with moduli $4,8,4*,8*,8*$ which are not sufficient to construct a covering, proving the claim.

Thus, ${\mathcal M_0'}$ contains at most one of the moduli $5*,10*,40*$. Allocating just one of $5*, 10*, 40*$ to ${\mathcal M_0'}$ does not help construct ${\mathcal C_0'}$. For example,  if $5*$ is in ${\mathcal M_0'}$ and $10*,40*$ are not in ${\mathcal M_0'}$, again, using Corollary \ref{cor:3} we can discard from ${\mathcal C_0'}$ the congruences with moduli $5,10,20,40,5*$. To complete ${\mathcal C_0'}$ we need to assign to ${\mathcal M_0'}$ one of the moduli $4*,8*,8*$. It is possible to contruct the covering ${\mathcal C_0'}$ by assigning to it one congruence modulo $8*$ and it is the efficient way to do it. (The coverings ${\mathcal C_1'}$ and ${\mathcal C_2'}$ cannot benefit from swapping  with  ${\mathcal C_0'}$ a congruence modulo $4$ with a congruence modulo $8$). So, one modulus $8*$ is allocated to ${\mathcal M_0'}$.

Next, we analyze how to split the moduli $5*,10*,40*$  among  ${\mathcal M_1'}$ and  ${\mathcal M_2'}$. 

We claim that both moduli $5*$ and $10*$ are in ${\mathcal M_2'}$. Assume ohterwise. Then, ${\mathcal M_2'}$ contains at most the moduli  $4,8,5,10,20,40, 4*,8*, 40*$ and one of the moduli $5*,10*$. Apply Lemma \ref{lem:1} to 
the congruences in the above list which are multiples of $5$. The only way to have $D$ nonempty set is to have $|1  | 2 | 4 \mbox{ and }8 | \mbox{ either } 1 \mbox { or } 2 | 8|$ in the five `bins' (here we use use $|$ as a separator between the bins and for brevity, instead of writing `the congruence modulo $2$ is in certain bin' we just write `$2$ is in the bin'). Also, note that in Case II, $4$ and $8$ are in the same `bin'  and we have either $1$ or $2$  in a `bin' depending on whether $5*$ or $10*$ is in ${\mathcal M_2'}$. Thus, 
we can replace the congruences with moduli which are multiples of $5$ by a single congruence modulo $8$. Now, we need to construct a covering with moduli  from the list $4,8,4*,8*,8$ which is impossible. 

Thus, we need to allocate both $5*$ and $10*$ to ${\mathcal M_2'}$. There are two subcases depending on how we allocate $40*$. 
\vskip 5pt
\noindent
{\bf Subcase A:} The modulus $40*$ is in ${\mathcal M_1'}$. 
\vskip 5pt
In this subcase ${\mathcal M_2'}$ contains the moduli  $4,8,5,10,20,40, 5*,10*$ and some of the the moduli $4*,8*$. We apply Lemma \ref{lem:1} again to the congruences with moduli $5,10,20,40, 5*,10*$. In this case $D$ is nonempty 
only if we have $|1 | 2 | 4 \mbox{ and } 8 | 1 | 2|$ in the 'bins' ( again, $4$ and $8$ have to be in the same `bin').  
Therefore, we can replace the congruences with moduli $5,10,20,40, 5*,10*$ by two congruences with moduli $4,8$ respectively. 
Now, we are left with moduli $4,8,4,8$ and some of $4*,8*$. So, we need to allocate $4*$ to ${\mathcal M_2'}$. This leaves for ${\mathcal M_1'}$ the moduli $4,8,5,10,20,40,3*,6*,12*,15*,8*,40*$.  We apply Lemma \ref{lem:1} again, this time to the congruences with moduli $5,10,20,40,15*,40*$. Again, $D$ is nonempty only if the content of the `bins' is $|1 | 2 | 4 \mbox{ and }8| 3 | 8|$ (in some order). Thus, we can replace the congruences with moduli $5,10,20,40,15*,40*$ by a single congruence modulo $24$.
We are left with moduli $4,8,3*,6*,12*,24*,8*$. We already proved in Case 1, Subcase B that it is not possible to construct a covering with moduli from the last list. The same proof works word by word here, too, so we are done with this subcase.
\vskip 5pt
\noindent
{\bf Subcase B:} The modulus $40*$ is in ${\mathcal M_2'}$. 
\vskip 5pt
In this subcase, ${\mathcal M_2'}$ contains the moduli  $4,8,5,10,20,40, 5*,10*,40*$ and some of the the moduli $4*,8*$. We apply Lemma \ref{lem:1} again to the congruences with moduli divisible by $5$. The congruences reduced modulo $5$ have moduli $1,2,4,8,1,2,8$ respectively. We need to place seven congruences in five 'bins', so at least one 'bin' will contain only one  congruence modulo $2, 4$, or $8$. Thus, one can replace the congruences with moduli $5,10,20,40, 5*,10*,40*$ by a single congruence modulo $2$. This leaves ${\mathcal M_2'}$ with moduli $2,4,8$ and some of $4*,8*$, so we allocate $8*$ to  ${\mathcal M_2'}$. Now, ${\mathcal M_1'} = [4,8,5,10,20,40,3*,6*,12*,15*,4*]$. We apply Corollary \ref{cor:3}  to the congruences with moduli $5,10,20,40,15*$.  Since the congruences with moduli $20$ and $40$ in ${\mathcal C'}$ are in the same class modulo $5$, we can discard the congruences  with moduli $5,10,20,40,15*$. This leaves us 
the moduli $4,8,3,6,12,4$. We apply Corollary \ref{cor:3} again to replace the congruences with moduli $3,6,12$ by a single congruence modulo $4$. Finally, we are left with moduli $4,8,4,4$ which is not sufficient to construct a covering. 

Having exhausted all cases, we obtain the proof of the theorem.
\end{proof}

\section{Covering systems with minimum least common multiple of the moduli} \label{sec:4}

In this section we solve the problem of minimizing the least common multiple of the moduli of a distinct covering system with a fixed minimum modulus $m$ in the cases $m=3$ and $m=4$. 

First, we need a lemma which will help us reduce the number of cases we need to consider. This lemma is Theorem 1 of Simpson and Zeilberger \cite{SimZeil} with the extra condition that the minimum modulus does not change.

\begin{lemma}  \label{lem:4}
Let ${\mathcal C}$ be a distinct covering with minimum modulus $m$, and least common multiple of the moduli  $L = L_1 q^{\alpha}$,  where $q$ is  a prime such that $q \nmid L_1$. Suppose $p$ is a prime which does not divide $L_1$, and $p$ is less than  $q$. Then, one can construct 
a distinct covering ${\mathcal C_1}$ with the same minimum modulus $m$, and least common multiple of the moduli  $L_1  p^\alpha$. 
\end{lemma}

\begin{proof} 

This proof relies on the `coordinate' notation for congruences we introduced. 

Let $q$ be the $j$th prime in the prime factorization of $L$. Write all congruences in ${\mathcal C}$ in `coordinate' notation. 
We  keep in ${\mathcal C_1}$ the congruences in which all base $q$ digits in the $j$th `coordinate'  are $\leq p-1$ with no change,  and discard the remaining congruences.
In the congruences which survived, we interpret the $j$th component modulo $p$. We claim that the  congruences in ${\mathcal C_1}$ form  a covering with least common multiple of the moduli  $L_1  p^\alpha$.

First, consider a residue class $r_1$ modulo  $L_1  p^\alpha$. Write $r_1$ in `coordinate' notation. Note that all digits in the $j$th position do not exceed $p-1$. The residue class $r_1$  corresponds to a residue class $r$ modulo  $L_1  q^\alpha$ where we have kept all digits in all positions the same (but interpreted the digits in $j$th position modulo $q$). 
Since ${\mathcal C}$ is a covering, there is a congruence $c$ in ${\mathcal C}$ which covers the residue class $r$. The congruence  $c$ corresponds to a congruence $c_1$ in ${\mathcal C_1}$, where both congruences have the same digits in all positions in `coordinate' notation.  
Clearly, $c_1$ covers $r_1$, so ${\mathcal C_1}$ is a covering. Moreover, ${\mathcal C}$ is a distinct covering, so by construction ${\mathcal C_1}$ is a distinct covering (all we do is replace $q$ by $p$ in the prime factorization of the moduli and discard some congruences). 

With the way our construction works, it is possible that the congruence with modulus $m$ in ${\mathcal C}$ may get discarded. However, note that if $x \equiv r_i$ (mod $n_i)$, $i = 1,\ldots ,k$ is a covering, then for any integer constant $a$, 
$x \equiv r_i+a$ (mod $n_i)$, $i = 1,\ldots ,k$ is also a covering since the set of all integers is  invariant to translation by an integer. Thus, if needed, we can translate the covering  ${\mathcal C}$ by a suitable constant $a$ (replace all $r_i$s by $r_i+a$), so that all digits   of the `coordinate' representation of the congruence  modulo $m$ in  ${\mathcal C}$, are zero (make the congruence $x \equiv 0$ (mod $m)$). In this way, the congruence modulo $m$ will not be discarded and there will be corresponding congruence modulo $m$ in ${\mathcal C_1}$.
\end{proof}
For example, consider the  covering ${\mathcal C}$ with $L = 80 = 2^4 5$, $(1)$, $(01)$, $(001)$, $(0001)$, 
$(*|\  4)$, $(0|\ 3)$, $(00|\ 2)$, $(000|\ 1)$, $(0000|\  0)$.

Proceeding as in the proof of Lemma \ref{lem:4}  with $q=5$ and $p=3$, we get the   covering ${\mathcal C_1}$ with $L = 48 = 2^4 3$, $(1)$, $(01)$, $(001)$, $(0001)$, $(00|\ 2)$, $(000|\  1)$,  $(0000|\  0)$. 

Now, we turn to Theorem \ref{lcm}. Erd\H{o}s constructed a covering ${\mathcal C}$ with least modulus $m = 3$. Krukenberg \cite{Krukenberg} also constructed a covering ${\mathcal C}$ with least modulus $m = 3$,  $L({\mathcal C}) = 120$,  without using the moduli $40$ and $120$. 
Here is a covering with the above properties $(11)$, $(101)$, $(*|\ 2 )$, $(0|\  1)$, $(100 |\ 1)$, $(10|\ 0)$,  $(*|\ *|\ 4)$, $(0|\ *|\ 3)$, $(*|\ 0|\ 2)$, $(0|\  0 |\  1)$, $(01|\ *|\ 0)$, $(00|\ 0|\ 0)$.
Next, we prove Theorem \ref{lcm}.

\begin{proof}

First we deal with part (i), the case $m=3$. 
We need to show that if $n$ is less than $120$ there is no distinct covering having as moduli only   divisors of $n$ which are at least $3$. Since the only $n$ less than  $120$ for which 
${\displaystyle \sum_{d|n, d \geq 3} \frac{1}{ d} \geq 1}$ are $24,36,48,60,72,84,90, 96$, and $108$, we only need to examine the numbers in this list. 

We can eliminate some cases using the work of Krukenberg on  coverings with least common multiple of the moduli of the form $2^a 3^b$, see Theorem  \ref{Krukenberg2}.
As proved by Krukenberg,  there is no covering with $m=3$, and $L= 24$, or $36$, or $48$, or $72$, or $96$, or $108$. What is left is to consider the cases when $m = 3$ and $L =60$, or $84$, or $90$.
By Lemma \ref{lem:4}, if there is a covering with $m = 3$ and $L = 84$, then there is a covering with $m = 3$ and $L = 60$. 

To finish the proof  we need to show that there is no covering with $m = 3$ and $L=60$ or $L=90$. 

First, assume that there is a covering ${\mathcal C}$  with $m = 3$ and $L=60$. Then, the moduli of the congruences are $4$; $3$, $6$, $12$; $5$, $10$, $20$; $15$, $30$, and $60$.  Reduce the covering modulo $3$. We have to construct three coverings with shared moduli: $4$, $5$, $10$, $20$, and moduli used by just one covering: $1*$, $2*$, $4*$, $5*$, $10*$, $20*$. Consider the covering, say ${\mathcal C_1}$  
which does not include $1*$, and includes at most one of $5*$, $10*$, and $20*$. It's moduli are 
$4$, $5$, $10$, $20$,  some  of $2*, 4*$, and at most one of  $5*, 10*, 20*$. We can discard from ${\mathcal C_1}$  all  congruences with moduli which are multiples of $5$ (there at most four of them). Thus, ${\mathcal C_1}$ includes both congruences with moduli $2*$ and $4*$. Let ${\mathcal C_2}$ be the covering including the congruence modulo 
$1*$. Then the moduli of the congruences in ${\mathcal C_3}$ ar at most $4$, $5$, $10$, $20$, $5*$, $10*$, $20*$.  The sum of the reciprocals of these moduli is at most $.95 < 1$, so a covering with with $m = 3$ and $L=60$ does not exist. 

Finally,  assume that there is a covering ${\mathcal C}$  with $m = 3$ and $L=90$. Then, the moduli of the congruences are $3$, $6$; $9$, $18$; $5$, $10$; $15$, $30$;  $45$, and $90$. Reduce the covering modulo $5$. We obtain five
coverings with shared moduli: $3$, $6$, $9$, $18$, and moduli used by just one covering: $1*$, $2*$, $3*$, $6*$, $9*$, $18*$. Consider the two  coverings that do not contain  any congruences modulo $1*$, $2*$, or $3*$. Since the sum of the reciprocals of $3, 6, 9, 18$ is $2/3$, both coverings need 
all three moduli  $6*$, $9*$, $18*$. Thus, a covering with $m = 3$ and $L=90$ does not exist completing the proof of part (i) of the theorem.

Now, we turn to part (ii), the case $m=4$. 

Krukenberg \cite{Krukenberg} constructed a covering ${\mathcal C}$ with $m = 4$ and $L({\mathcal C}) = 360$.  Here is a covering which uses as moduli all divisors of $360$ which are at least  $4$, except $360$. It is $(11)$, $(101)$, $(0|\ 2)$, 
$(100|\ 2)$, $(01|\ 1)$, $(*|\ 02)$, $(0|\ 01)$, $(100|\ 01)$, $(10|\ 00)$, $(*|\ *|\ 4)$, $(0|\ *|\ 3)$, $(100|\ *|\ 3)$, $(00|\ *|\ 2)$, $(*|\ 1|\ 0)$, $(0|\ 1|\ 1)$, $(10|\ 1|\ 1)$, $(100|\ 1|\ 2)$, $(*|\ 00|\ 0)$, $(0|\ 00|\ 1)$, $(01|\ 00|\ 2)$.

We need to show that if $n$ is less than $360$ there is no covering using only distinct divisors of $n$ which are at least $4$. Since the only positive integers $n$ less than  $360$ for which 
${\displaystyle \sum_{d|n, d \geq 4} \frac{1}{ d} \geq 1}$ are $120  ,168 , 180 ,240 ,252,  280 ,288, 300$, and $336$ we only need to examine these values of $n$. 

Since, $120 | 240$, it is sufficient to show that $240$ does not work. 

Using Lemma \ref{lem:4}  we can reduce the cases $n=168= 2^3 \cdot 3 \cdot 7$, $n=252=2^2 \cdot 3^2 \cdot 7$, and  $n = 336 = 2^4 \cdot 3 \cdot 7$ to $n=120$, $n=180$,  and $n=240$ respectively. 

As proved by Krukenberg,  Theorem \ref{Krukenberg2}, $n = 288 = 2^5 3^2$ does not work either. 

Let us consider the case $n=280$. The sum of the reciprocals of the divisors of $280$ which are at least $4$ is $1.0714\ldots$. However the congruences modulo $4$, $5$, and $7$ cover a portion of the integers with density 
$1 - \frac{3}{4} \cdot  \frac{4}{5} \cdot \frac{6}{7}= \frac{17}{35}$. Since $ \frac{1}{4}  + \frac{1}{5}   + \frac{1}{7}  -  \frac{17}{35} = .1071\ldots$, there is no covering with $m = 4$ and $L = 280$.  

Next, let $n =300$. The sum of the reciprocals of the divisors of $300$ which are at least $4$ is $1.06$. However the intersection of the congruences modulo $4$, $5$, and $15$ is at least $ \frac{1}{15} = .0666\ldots $. Therefore, 
there is no covering with $m = 4$ and $L = 280$. 

We are left with two remaining cases: $n=180$ and $n=240$. We consider each case in a separate lemma. 

\begin{lemma} \label{lem:5}
There is no distinct covering with $m = 4$ and $L=180$. 
\end{lemma}

\begin{proof}
Assume that ${\mathcal C}$ is a covering with moduli $4$;  $6$, $12$; $9$, $18$, $36$; $5$, $10$, $20$; $15$, $30$, $60$; $45$, $90$, and $180$. Let $S$ be the set of congruences with moduli $4$, $6$, $12$, $9$, $18$, $36$. 
Note that the density of the integers covered by congruences in $S$ is at most $\frac{2}{3}$. Indeed,  the sum of the reciprocals of the moduli of congruences in $S$  is $\frac{25}{36}$ and  the set of integers covered by the congruences modulo $4$ and modulo $9$ intersect. 

Next, reduce ${\mathcal C}$ modulo $5$. 
We need to construct five coverings with common moduli $4$, $6$, $12$, $9$, $18$, $36$  and moduli used by just one covering: $1*$, $2*$, $4*$, $3*$, $6*$, $12*$, $9*$, $18*$, $36*$.

Consider the three coverings containing the congruences with moduli $1*$,  $2*$ , $3*$. One can see that either the congruence modulo $2$ and $S$ do  not form a covering or 
the congruence modulo $3$ and $S$ do  not form a covering. Otherwise, the set uncovered by $S$  is inside a residue class modulo $2$ and inside a residue class modulo $3$, that is, inside a residue class modulo $6$. This is not possible since the set uncovered by $S$ has density at least $1/3$. 

Thus, the three coverings containing the congruences with moduli  $1*$,  $2*$ , $3*$ contain at least one more congruence. We can assume it has modulus $36*$ (all other $*$ moduli are divisors of $36$). So, the fourth covering and the fifth covering need to split the moduli
$4*$, $6*$, $12*$, $9*$, $18*$. Now, $\frac{1}{4} +\frac{1}{6}    + \frac{1}{12}  + \frac{1}{9}  + \frac{1}{18}  = \frac{2}{3} $. Recall that the set uncovered by $S$ has density at least $1/3$. Thus,  the only possible way to construct the remaining two  coverings is  if one covering uses moduli $4*$ and $12*$ and the other covering uses  $6*$, $9*$, $18*$. Therefore,  we need to be able to construct a covering using the moduli in the list $[4, 4,6, 12, 12,9, 18, 36]$.  By Corollary \ref{cor:3}  we can replace the congruences with moduli  $9$, $18$, $36$ by  a single congruence modulo $12$. The moduli of the congruences of  resulting covering are in the list $[4, 4,6, 12, 12,12]$. 
Since the sum of the reciprocals of the elements of the  list $[4, 4,6, 12, 12,12]$  is less that one,  it is not possible to construct at least one of the five coverings we needed to construct. 
\end{proof}

The proof of the next lemma is somewhat complicated. However, we expect that the methods used in the proof of the lemma will be useful when analyzing coverings with least common multiple of the moduli of the form $2^a 3^b 5^c$. 
\begin{lemma} \label{lem:6}
 There is no distinct covering with $m = 4$ and $L=240$.
\end{lemma}
 
\begin{proof}

In the proof of this lemma we will use the notation $(n,r_n)$ to denote the congruence $x \equiv r_n\ (\mbox{mod} \ n)$.

  Assume  that there is a covering\\ ${\mathcal C} = \{ (n,  r_n)\ |\ n \in \{ 4,8,16, 6,12,24,48,5,10,20,40,80,15,30,60,120,240\} \}$. We introduce notation for some of the parts of ${\mathcal C}$. 
Let ${\mathcal C_2} = \{  (n,  r_n)\ |\ n \in \{ 4,8,16 \} \}$, \\${\mathcal C_3} = \{  (n,  r_n)\ |\ n \in \{ 6,12,24,48 \} \}$, ${\mathcal C_5} = \{  (n,  r_n)\ |\ n \in \{ 10,20,40,80 \} \}$, and\\ 
${\mathcal C_{15}} = \{  (n,  r_n)\ |\ n \in \{ 15,30,60,120,240 \} \}$.\\
\noindent Also, let $R$ be the set  of residue classes modulo $16$ which are  not covered by the congruences in ${\mathcal C_2}$. Thus, $R$ is a set of $9$ residue classes modulo $16$ represented by $9$ integers between $0$ and $15$. 
Let $R_0 = R \cap \{ x \equiv 0\ (\mbox{mod} \ 2) \}$ and 
  $R_1 = R \cap \{ x \equiv 1\ (\mbox{mod} \ 2) \}$.

For each $r \in R$ denote by $a_3(r)$ the number of residue classes modulo $48$ of the form\\
\noindent  $ x \equiv r \ (\mbox{mod} \ 16) , \  x \equiv a \ (\mbox{mod} \ 3)$, which are covered by ${\mathcal C_3}$. One way to visualize this is that the residue class 
$(r  \ \mbox{mod} \ 16)$ splits into three fibers modulo $48$. The quantity $a_3(r)$ counts how many of these fibers are covered by  ${\mathcal C_3}$.

Similarly, for each $r \in R$ denote by $a_5(r)$ the number of residue classes modulo $80$ of the form
\noindent $ x \equiv r \ (\mbox{mod} \ 16) , \  x \equiv b \ (\mbox{mod} \ 5)$, which are covered by ${\mathcal C_5}$. 

Then,  the number of residue classes modulo $240$ which are not covered by any of  the congruences in  ${\mathcal C_2}$,  ${\mathcal C_3}$,  ${\mathcal C_5}$, nor by the congruence $(5, r_5)$ is 

\begin{equation} \label{left mod 15}
L:= \sum_{r \in R} (3 - a_3(r))(4-a_5(r))
\end{equation}

Note that the congruences in ${\mathcal C_{15}}$ can cover at most $5$ residue classes modulo $240$ which are in the residue class $(r  \ \mbox{mod} \ 16)$. Thus, for each $r \in R$ we have 
$(3-a_3(r))(4 - a_5(r)) \leq 5$. Clearly, $(3-a_3(r))(4 - a_5(r)) \neq 5$.

\noindent So, for each $r \in R$ we have

\begin{equation} \label{at most 4}
 (3 - a_3(r))(4-a_5(r)) \leq 4.
\end{equation}

Furthermore, for each $r \in R$,

\begin{equation} \label{2*3>5}
\mbox{ if } \quad a_3(r)a_5(r) \neq 0, \mbox{ then } a_3(r)a_5(r) \geq 2. 
\end{equation}

One more observation: Suppose $r_1 \in R$, $r_2 \in R$, and  $r_1 \not\equiv r_2 \ (\mbox{mod} \ 2)$. Then the number of residue classes modulo $240$ which are either
$\equiv r_1  \ (\mbox{mod} \ 16)$ or $\equiv r_2  \ (\mbox{mod} \ 16)$ and can be covered by ${\mathcal C_{15}}$ is at most $6$. Indeed, the congruence modulo $15$  can cover at most two such classes, and each of 
$(30,r_{30})$, $(60, r_{60})$, $(120, r_{120})$, and $(240, r_{240})$ can cover at most one. So, in this case 

\begin{equation} \label{odd even}
 (3 - a_3(r_1))(4-a_5(r_1))+  (3 - a_3(r_2))(4-a_5(r_2)) \leq 6.
\end{equation}

We can rewrite \eqref{left mod 15} as 

\begin{equation} \label{S3 S5}
 L = \sum_{r \in R} (12 - 4a_3(r) - 3a_5(r)+a_3(r)a_5(r)) = 108 - 4S_3 - 3S_5 + O,
\end{equation}

where ${\displaystyle S_3 = \sum_{r \in R} a_3(r)}$,  ${\displaystyle S_5 = \sum_{r \in R} a_5(r)}$, and  ${\displaystyle O = \sum_{r \in R} a_3(r)a_5(r)}$. 

The quantity $O$ measures the amount of overlap between ${\mathcal C_3}$ and ${\mathcal C_5}$. Ideally,  we want $O$ to be small (cover one set of $r$'s by ${\mathcal C_3}$ and a different set of $r$'s  by ${\mathcal C_5}$), while $S_3$ and $S_5$ are large, that is,  
cover a lot without much overlap. Al least in the case of this lemma, this proves impossible.

Next, we get bounds for $S_3$ and $S_5$. 

For $n \in \{ 2,4,6,8,16 \}$ define ${\displaystyle M_n = \max_{0 \leq j < n}  |\{ R \cap \{  x \equiv j\ (\mbox{mod} \ n) \} \}|}$. Here  $M_n$ is the size of the largest portion of $R$ in a residue class modulo $n$. 

Then, the congruence $(6, r_6)$ can contribute at most $M_2$ to $S_3$, the congruence $(12, r_{12})$ - at most $M_4$, etc. 

Thus, $S_3 \leq M_2 + M_4 + M_8 + M_{16}$. Similarly,  $S_5 \leq M_2 + M_4 + M_8 + M_{16}$.
Define $D_3 = (M_2 + M_4 + M_8 + M_{16}) - S_3$ and   $D_5 = (M_2 + M_4 + M_8 + M_{16}) - S_5$.  

In certain sense, $D_3$ and  $D_5$  measure the difference between the largest amount we could possibly cover,  and what we cover in reality with 
${\mathcal C_3}$ and  ${\mathcal C_5}$,  respectively. For example, if $R$ consists of $1$ class $r$ such that $ r  \equiv 0\ (\mbox{mod} \ 2)$,  and $8$ classes $r_1$ such that $ r_1 \equiv 1\ (\mbox{mod} \ 2)$, then if we have a congruence $(6, r_6)$ with 
	$r_6 \equiv 0\ (\mbox{mod} \ 2)$, then $D_3 \geq 7$ (we could have covered $8$ residue classes and covered just $1$ instead.)

Also, the number of residue classes modulo $240$ which  can be covered by ${\mathcal C_{15}}$ does not exceed $9  +M_2 + M_4 + M_8 + M_{16}$.
Therefore, if ${\mathcal C}$ is a covering, then  $ L  \leq  9  +M_2 + M_4 + M_8 + M_{16}$. Recall that $L$ is the number of residue classes modulo $240$ which are not covered by any of  the congruences in  ${\mathcal C_2}$,  ${\mathcal C_3}$,  ${\mathcal C_5}$, nor by the congruence $(5, r_5)$. These classes need to be covered by  ${\mathcal C_{15}}$.

\noindent Define  $D_{15} = 9+ (M_2 + M_4 + M_8 + M_{16})- L$. Since by assumption ${\mathcal C}$ is a covering, $D_{15} \geq 0$.

Using \eqref{S3 S5}, we get 
$9 + 8(M_2 + M_4 + M_8 + M_{16}) \geq 108 + 4D_3 + 3D_5 + D_{15} + O.$

Since, $M_4 \leq 4$, $M_8 \leq 2$, and $M_{16} \leq 1$, we obtain

\begin{equation} \label{M2}
8M_2 \geq 43  + 4D_3 + 3D_5 + D_{15} + O.
\end{equation}

Next, we consider several cases, depending on the structure of ${\mathcal C_2}$.Without loss of generality  we can assume $r_4 \equiv 0 \pmod 2$. As we noted above, if $\{ (n, r_n) | n \in L \}$, where $L$ is a list of moduli, is a covering, then for any integer $a$, 
 $\{ (n, r_n+a) | n \in L\}$ is also a covering. 

\vspace{0.10 in} \noindent
{\bf Case I.} $r_8 \equiv  r_{16}   \equiv 1\ (\mbox{mod} \ 2)$. 

In this case, $|R_0|=4$, $|R_1|=5$, and $M_2 = 5$. 

From \eqref{M2}  we get $0 \geq 3 +  4D_3 + 3D_5 + D_{15} + O$. Since $D_3 \geq 0$, $D_5 \geq 0$, $D_{15} \geq 0$, and $O \geq 0$, we get a contradiction. 
There is no covering in Case I. 

\vspace{0.10 in} \noindent
{\bf Case II.} $r_8   \equiv 1 \ (\mbox{mod} \ 2)$ and $ r_{16}  \equiv 0 \ (\mbox{mod} \ 2)$.

In this case, $|R_0|=3$, $|R_1|=6$, and $M_2 = 6$. 

From \eqref{M2} we get $5 \geq   4D_3 + 3D_5 + D_{15} + O$. Thus, $D_3 \leq 1$ and $D_5 \leq 1$. Hence, $r_6  \equiv 1 \ (\mbox{mod} \ 2)$ and $r_{10}  \equiv 1 \ (\mbox{mod} \ 2)$. We obtain that $a_3(r) \geq 1$ and $a_5(r) \geq 1$ for all $r \in R_1$, 
so $O \geq 6$, a contradiction in this case, too.

\vspace{0.10 in} \noindent
{\bf Case III.} $r_8  \equiv 0 \ (\mbox{mod} \ 2)$ and $ r_{16}  \equiv 1 \ (\mbox{mod} \ 2)$.

In this case, $|R_0|=2$, $|R_1|=7$, and $M_2 =7$. 

From \eqref{M2}  we get $13 \geq   4D_3 + 3D_5 + D_{15} + O$. Therefore $D_3 \leq 3$ and $D_5 \leq 4$. Again, $r_6 \equiv r_{10} \equiv 1 \pmod 2$. So, 
$a_3(r) \geq 1$ and $a_5(r) \geq 1$ for all $R \in R_1$. By \eqref{2*3>5}, $a_3(r)a_5(r) \geq 2$ for all $r \in R_1$. Therefore, $O \geq 14$, so a covering does not exist in this case, too.

\vspace{0.10 in} \noindent
{\bf Case IV.} $r_8 \equiv  r_{16} \equiv 0 \pmod 2$.

Here, $|R_0|=1$, $|R_1|=8$, and $M_2 =8$. So, $R_0 = \{ r_0 \  \mbox{(mod 16)} \}$ where $r_0$ is an integer in 
 $\{ 0,2,4,6 \}$.  

In this case, we can cover a lot with ${\mathcal C_3}$ and ${\mathcal C_5}$ but the overlap between them is too big and again, we fall short of constructing a covering. 

First, note that ${\displaystyle  \sum_{r \in R_1} a_3(r) \leq 8 + 4 + 2 + 1 = 15}$. Therefore, there exists $r_1 \in R_1$ 
such that $a_3(r_1) \leq 1$. By \eqref{at most 4}, we get $a_5(r_1) \geq 2$. Thus, $r_{10} \equiv r_{20} \equiv 1 \pmod 2$ (if any of the congruences in ${\mathcal C_5}$ are used to cover $R_0$, they should be the ones with  the largest moduli since $|R_0| = 1$). 

Since, $r_{10} \equiv r_{20}  \equiv 1 \ (\mbox{mod} \ 2)$, we have $a_5 (r_0) \leq 2$, and by  \eqref{at most 4}  we get $a_3(r_0) \geq 1$. Therefore, 
$r_{48}  \equiv 0 \ (\mbox{mod} \ 2)$. 

Similarly as above, ${\displaystyle  \sum_{r \in R_1} a_5(r) \leq  15}$. Therefore, there exists $r_1' \in R_1$ such that $a_5(r_1') \leq 1$. By \eqref{at most 4},  we get $a_3(r_1') \geq 2$. Thus, $r_{6} \equiv r_{12}  \equiv 1 \ (\mbox{mod} \ 2)$. So, 
$a_3 (r_0) \leq 2$. We proved above that $a_3(r_0) \geq 1$, so $a_3 (r_0)$  is either $1$ or $2$.

Assume that $a_3 (r_0)=1$. Then \eqref{at most 4}  implies $a_5(r_0) \geq 2$, so $r_{40} \equiv r_{80} \equiv 0 \pmod 5$. Hence,  
$a_5 (r) \leq 2$ for all $r \in R_1$. This implies $(3 - a_3(r_1))(4 - a_5(r_1)) \geq 4$. Also, 
$(3 - a_3(r_0))(4 - a_5(r_0)) = 4$. Thus, $(3 - a_3(r_1))(4 - a_5(r_1)) + (3 - a_3(r_0))(4 - a_5(r_0)) \geq 8$, which contradicts \eqref{odd even}. 

So, $a_3(r_0)=2$, and $r_{24} \equiv r_{48}  \equiv 0 \ (\mbox{mod} \ 2)$. 

 Next, let $R_1' = \{ r \in R_1 \ | \  r \not\equiv r_{12} \pmod 4 \}$. For all $r \in R_1'$ we have $a_3(r) = 1$.

\noindent Also,   ${\displaystyle  \sum_{r \in R_1'} a_5(r) \leq   4 + 4 + 2 + 1 = 11}$, so there exists $r_1^* \in R_1'$ 
with $a_5(r_1^*) \leq 2$.\\
\noindent  Then $(3 - a_3(r_1^*))(4 - a_5(r_1^*))  \geq 4$. By \eqref{odd even} we get $( (3 - a_3(r_0))(4 - a_5(r_0))  \leq 2$. Thus, 
$a_5 (r_0) = 2$, and $r_{40} \equiv r_{80} \equiv 0 \pmod 2$. 

We have allocated all congruences in ${\mathcal C_3}$ and ${\mathcal C_5}$ to $R_0$ and $R_1$ (both $R_0$ and $R_1$ get two congruences from ${\mathcal C_3}$ and two 
from ${\mathcal C_5}$). 

Since $M_2 = 8$, \eqref{M2}  becomes 

\begin{equation} \label{21}
21 \geq 4D_3 + 3D_5 + D_{15} + O.
\end{equation} 

However, $D_3 \geq 1$, since $(24, r_{24})$ covers just one class modulo $48$, and  $D_5 \geq 1$ since  we did not use 
$(40, r_{40})$ in the most efficient way either. 

Also, $a_3(r) a_5(r) \neq 0$ for all $r \in R$, so by \eqref{2*3>5}  $a_3(r) a_5(r) \geq  2$ for all $r \in R$, and $O \geq 18$. 
Substituting in (13) we get $21 \geq 4 + 3 + 18$, a contradiction. 
\end{proof}

We just considered the last remaining case and
this concludes the proof of Theorem \ref{lcm}.
\end{proof}

\section{Open problems and further work} \label{open pr}

Recall that Krukenberg constructed a distinct   covering system with least modulus $5$ and largest modulus $108$. He also conjectured that one can not replace $108$ by a smaller constant. 

\vspace{0.1 in} 
\noindent {\bf Problem 1.} Prove or disprove that if the least modulus of a distinct covering system is $5$, then its largest modulus is at least $108$. 

\vspace{0.1 in} \noindent
We can show that if the least modulus of a distinct covering system is $5$, then its largest modulus is at least $84$. However, the result is too weak and the proof too long, to be included in this paper.

\vspace{5 pt}
Krukenberg also provided a description of the covering systems with least common multiple of the moduli of the form $2^a 3^b$, see Theorem \ref{Krukenberg2}. 

\vspace{0.1 in} \noindent
{\bf Problem 2.} Describe the distinct covering systems with least common multiple of the moduli of the form $2^a 3^b 5^c$ where $a$, $b$, and $c$ are integers such that $a \geq b \geq c \geq 1$. 

\vspace{0.1 in} \noindent
Using Krukenburg's results and the results of this paper one can find such a description when $m=2,3,4$ with one exception. Extra work is needed to show that there is no distinct covering system with $m=4$ and $L=900$. 
The open case is when $m \geq 5$. 

\vspace{5 pt}
Furthermore, Krukenberg constructed a distinct  covering system with $m=5$ and\\ $L=1440$. 

\vspace{0.1 in} \noindent
 {\bf Problem 3.} Prove or disprove that if the least modulus of a distinct covering system is $5$, then the least common multiple of its moduli is at least $1440$.

\vspace{5 pt}
Krukenberg also constructed a distinct covering system not using the modulus $3$, with all moduli  squarefree integers. It is not known whether there exists a  distinct covering system  with squarefree moduli and least modulus $3$. 

\vspace{0.1 in} \noindent
 {\bf Problem 4.} Prove or disprove that the least modulus of  any  distinct covering system  with squarefree moduli is $2$. 

\vspace{0.1 in} \noindent
Solving Problem 4 will give us a complete solution of the {\it minimum modulus problem} in the squarefree case. 

\vspace{0.1 in} \noindent
{\bf Problem 5.}  Find the largest integer $c$ such that there exists a finite set of congruences with distinct moduli with the property that every integer satisfies at least $c$ of the congruences. 
In other words, what is the largest number of times we can cover the integers by a finite system of congruences with distinct positive moduli?

\vspace{0.1 in} \noindent
For a positive integer $n$ let $c(n)$ be the largest number of times we can cover all integers using congruences with moduli $1,2,\ldots,n$ respectively. Clearly, $c(1)=1$. 
Also, $c(n) \leq c(n+1)$ for all positive integers $n$ (having more congruences allows us to cover more).  Furthermore, $c(n+1) \leq c(n) + 1$ for all $n$. Indeed, if certain integer is covered $c(n)$ times by 
certain congruences with moduli $1,\ldots,n$, it can be covered at most once more  by a congruence with modulus $n+1$.

Recall that by Theorem \ref{Krukenberg1} there is no distinct covering with moduli in the interval $[2,11]$, and there is a distinct covering with moduli $2,3,4,6,12$. Therefore, $c(2) = \cdots = c(11) = 1$, and $c(12) = 2$. 
Moreover, Krukenberg constructed a distinct  covering with least modulus $13$ and largest modulus $52562109600$. Therefore,  $c(52562109600) \geq 3$. 

Moreover, the sequence $\{ c_n \}_{n=1}^\infty $ is bounded. Recall that Balister et. al. \cite{Bal0} showed that the least modulus of any distinct covering system does not exceed $606000$. If we consider a system of congruences with moduli $1,2, \ldots , n$ respectively, where $n > 606000$, there will be an integer $m$ which is not covered by any of the congruences with moduli $606001, 606002, \ldots, n$. Even if $m$ is covered by each of the congruences with moduli $1,\ldots, 606000$, then $m$ will be covered $606000$ times. Thus, $c(n) \leq 606000$ for all $n$. 

Thus, ${\displaystyle c = \lim_{n \to \infty} c(n)}$ and Problem 5 is to find $c$. 

We have a heuristic based on several assumptions showing that $c$ is either $4$ or $5$.   So far, we know $3 \leq c \leq 606000$. 

\vspace{0.15 in}
\noindent
{\bf Acknowledgements} The authors greatly appreciate the help of Michael Filaseta whose numerous suggestions  helped improve the paper.

\vspace{0.15 in}
\noindent
{\bf Appendix 1}
\vskip 5pt
Here we provide some details on how we showed that $T_m < 1$ for all $m$ in $[51,606000]$. Using the methods below, we only needed to calculate $20$ values of $T_m$ in the interval $[69, 606000]$. 

Recall that ${\displaystyle T_m = \sum_{\substack{ m \leq n \leq 8m \\ P(n) < \sqrt{7m+1}}} \frac{1}{n}}$, where $P(n)$ denotes the largest prime divisor of $n$. 

First, we computed and stored $P(n)$ for all $n$ from $2$ to $8\cdot 606000=4848000$. 

Next note that $T_{m-1} \leq T_m + a_{m-1}$, where we define $a_{m-1}$ to be $\frac{1}{m-1}$  when $P(m-1) < \sqrt{7m-6}$, and we define $a_{m-1}$ to be $0$ when $P(m-1) \geq \sqrt{7m-6}$.

Indeed, \begin{align*}
	T_{m-1}  &= \sum_{\substack{ m-1 \leq n \leq 8m-8 \\ P(n) < \sqrt{7m-6}}} \frac{1}{n}  \\
	& = a_{m-1}  +  \sum_{\substack{ m \leq n \leq 8m-8 \\ P(n) < \sqrt{7m-6}}} \frac{1}{n}\\
	& \leq a_{m-1}  +  \sum_{\substack{ m \leq n \leq 8m \\ P(n) < \sqrt{7m+1}}} \frac{1}{n}. 
\end{align*}

So, we computed $T_{606000} = 0.6889312518\ldots$, and then using the inequality $T_{m-1} \leq T_m + a_{m-1}$ we get that $T_{605999} \leq T_{606000} +\frac {1}{605999}$. 
Iterating this method, we backtracked down to the last value of $m$ where the sum is less than $1$, which got us to $281678$. 
Thus, because $$\sum_{m=281678}^{605999} a_m + T_{606000}=0.9999971415\ldots<1,$$ we get that $T_m < 1$ for all $m \in [281678, 606000]$. 
Next we computed $$T_{281677}=0.6904637039\ldots,$$ and backtracked again; we got $$\sum_{m=133009}^{281677} a_m + T_{281678}=0.9999959662\ldots<1.$$
Through these jumps, we confirmed that $T_m <1$ for each $m$ in the intervals 
\begin{align*}
[281678, 606000], [133009, 281677], [64358, 133008], [3&2259, 64357], [16645, 32258], \\
[8721, 16644], [4775, 8720], [2731, 4774], [1649, 2730]&, [1028, 1648], [661, 1027], \\
[426, 660], [289, 425], [208, 288], [148, 207], [117, 147]&, [93, 116], [76, 92], [69, 75]\end{align*}

We had to compute $20$ values of $T_m$ to get to $m=69$. 
We then computed directly all values of $T_m$ for $m \in [3,68]$.

\vspace{0.15 in}
\noindent
{\bf Appendix 2}
\vskip 5pt

Here we provide an example of a distinct covering system with a congruence modulo $180$ and the remaining  moduli in $[4,56]$. 
The moduli we use are
$$4,8,16;6,12,24,48;9,18,36;5,10,20,40;
15,30;45,180;7,14,21,28,35,42,56,$$ 
where the semicolumns are used to separate the moduli involved in different stages of our argument below. 

The congruences modulo $4,8,16$ which we use are $(11)$, $(101)$, and $(1001)$. The uncovered set after the first  stage consists of a residue class modulo $2$, $(0)$, and a residue class modulo $16$, $(1000)$. 

Splitting modulo $3$, the uncovered set is $(0 |\ 0,1,2)$ and $(1000 |\ 0,1,2)$. 

Next, we use the congruences modulo $12,24,48$ to cover $(1000)$. The congruences modulo $12,24,48$ given by  $(10 | \ 0)$, $(100 |\ 1)$, and $(1000 |\ 2)$ accomplish this. 

We use the congruence modulo $6$ given by $(0 |\ 2)$. After the second stage, the uncovered set is $( 0 |\ 0,1)$. 

We use the congruences modulo $9,18,36$ to attack the residue class $(0 |\ 1)$, which is the same as $(0 |\ 10,11,12)$. We take the  congruences modulo $9,18,36$ to be $( * | 12)$, $(0 |\ 11)$, and $(01 |\ 10)$. 
The uncovered set after the third stage is $(0 |\ 0)$ and $(00 |\ 10)$. 

We split the uncovered set modulo $5$, to get 
$(0 |\ 0 |\ 0,1,2,3,4 )$ and $(00 |\ 10 |\ 0,1,2,3,4 )$. The congruences modulo $5$, $10$, and $20$ are $(* |\ * |\  4)$, $(0 |\ * |\ 3)$, and $(00 |\ * |\ 2)$.
Now, the uncovered set is $(0 |\ 0 |\ 0,1)$,  $(01 |\ 0 |\ 2 )$, and $(00 |\ 10 |\ 0,1 )$.  The congruence modulo $40$ is $(011 |\ * |\ 2 )$.
The congruences modulo $15$ and $30$ are $(* |\ 0 |\ 1)$ and $(0 |\ 0 |\ 0)$, and they cover $(0 |\ 0 |\ 0,1)$.  We are left with the uncovered set 
$(010 |\ 0 |\ 2 )$ and $(00 |\ 10 |\ 0,1 )$. We use the congruences modulo $45$ and $180$, $(* |\ 10 |\ 1 )$ and $(00 |\ 10 |\ 0 )$ to cover $(00 |\ 10 |\ 0,1 )$.
We are left with the single uncovered residue class $(010 |\ 0 |\ 2 )$ which we cover with the last seven congruences $(* |\ *|\ * |\ 6)$, $(0 |\ *|\ * |\ 5)$, $(* |\ 0 |\ * |\ 4)$, $(01 |\ *|\ * |\ 3)$, $(* |\ *|\ 2 |\ 2)$,
$(0 |\ 0|\ * |\ 1)$, and $(010 |\ *|\ * |\ 0)$.


\begin{thebibliography}{10}

\bibitem{Bal0}
Paul~Balister, B\'ela~ Bollob\'as,   Robert~Morris, Julian~Sahasrabudhe, Marius~Tiba, 
On the Erd\H{o}s covering problem: the density of the the uncovered set, preprint, 2018
 
\newblock Available at \url{https://arXiv.org/pdf/1811.03547.pdf}.

\bibitem{Bal1}
Paul~Balister, B\'ela~ Bollob\'as,   Robert~Morris, Julian~Sahasrabudhe, Marius~Tiba, The Erd\H{o}s–Selfridge problem with square-free moduli, \emph{Algebra Number Theory} \textbf{ 15} (2021),  609–-626.


\bibitem{Choi}
S.~L.~G.~Choi, Covering the set of integers by congruence classes of distinct moduli, \emph{Math.
Comp.} \textbf{ 25} (1971), 885-–895. 

\bibitem{ChurchHouse}
R.~ F.~Churchhouse, Covering sets and systems of congruences, {\it Computers in Mathematical
Research}, North-Holland,  1968,  20–-36.


\bibitem{Erd}
P.~Erd\H{o}s, On integers of the form $2^k+p$ and some related problems, \emph{Summa Brasil. Math. }
  \textbf{2}  (1950), 113–123.

\bibitem{ErdGra}
P. ~Erd\H{o}s and R.~L.~Graham, Old and new problems and results in combinatorial number theory, \emph{ Monographies de L’Enseignement Math\'ematique}, \textbf{28}, 1980.

\bibitem{Erdodd}
P. ~Erd\H{o}s, R\'esultats et probl\`emes en th\'eorie des nombres, \emph{S\'eminaire Delange-Pisot-Poitou (14e ann\'ee:
1972/73), Th\'eorie des nombres}, Fasc. 2, Exp. No. 24 , 7 pp., Secr\'etariat Math\'ematique, Paris, 1973.


\bibitem{Fil}
M.~Filaseta, K.~Ford, S.~Konyagin, C.~Pomerance, and G.~Yu, Sieving by large integers and covering
systems of congruences, \emph{J. Amer. Math. Soc.}, \textbf{20} (2007), 495-–517

\bibitem{Fil2}
M.~Filaseta, K.~ Ford, and S.~Konyagin, On an irreducibility theorem of A.~Schinzel associated with
coverings of the integers, \emph{Illinois J. Math.} \textbf{ 44} (2000), 633-–643.

\bibitem{Fil3}
M.~Filaseta, O.~Trifonov, and G.~Yu, Distinct covering systems, \emph{research seminar} (2006)

\bibitem{Gibson}
D.~J.~Gibson, A covering system with least modulus $25$, \emph{Math. Comp.} \textbf{ 78}  (2009),  1127–-1146. 

\bibitem{Hough}
R.~Hough, Solution of the minimum modulus problem for covering systems, \emph{Ann. Math.}, \textbf{181} (2015),
361-–382.

\bibitem{Krukenberg}
C.~E.~Krukenberg, Covering sets of the integers, \emph{Ph.D. thesis}, University of Illinois, Urbana-Champaign, 1971.

\bibitem{Morikawa}
Ryozo~Morikawa, On a method to construct covering sets, \emph{Bull. Fac. Liberal Arts Nagasaki
Univ.}  \textbf{22} (1981),  1–-11.

\bibitem{Nielsen}
P.~Nielsen, A covering system whose smallest modulus is $40$. \emph{J. Number Theory} \textbf{ 129},  (2009),  640-–666.

\bibitem{Owens}
Tyler~Owens, A Covering System with Minimum Modulus 42, \emph {Theses and Dissertations} \textbf{ 4329},   2014.

\newblock Available at \url{https://scholarsarchive.byu.edu/etd/4329}.

\bibitem{SimZeil}
R.~J.~Simpson and D.~Zeilberger, Necessary conditions for distinct covering systems with square-free 
moduli, \emph{Acta Arith.}\textbf{ 59} (1991), 59–-70.

\bibitem{Swift}
J.~D.~Swift, Sets of covering congruences, \emph{Bull. Amer. Math. Soc.} \textbf{60} (1954), 390.


\end{thebibliography}
\end{document}